\documentclass[11pt,english,oneside]{amsart}

\usepackage[breaklinks,colorlinks,linkcolor=blue,citecolor=blue,urlcolor=blue]{hyperref}
\usepackage[T1]{fontenc}
\allowdisplaybreaks[4]
\usepackage[latin9]{inputenc}
\usepackage{geometry}
\usepackage{indentfirst}
\usepackage{amssymb}
\usepackage{color}

\usepackage{booktabs}
\usepackage{array, caption, threeparttable}
\usepackage[font=small,labelfont=bf,labelsep=none]{caption}
\usepackage{graphicx,color}
\usepackage{amsmath, amssymb, graphics}

\usepackage{graphicx}

\makeatletter
\numberwithin{equation}{section} 
\numberwithin{figure}{section} 
  \@ifundefined{theoremstyle}{\usepackage{amsthm}}{}
  \theoremstyle{plain}
  
  \theoremstyle{plain}
  
  \theoremstyle{plain}
  
  \theoremstyle{remark}
  
  \theoremstyle{remark}
  
  \theoremstyle{plain}
  
\usepackage{geometry}

\def\exp{\hbox{\rm exp}}
\smallskip
\def\<{{\langle }}
\def\>{{\rangle }}

\makeatother

\usepackage{babel}

\def\exp{\hbox{\rm exp}}
\smallskip
\def\<{{\langle }}
\def\>{{\rangle }}

\makeatother

\usepackage[mathscr]{eucal}
\usepackage{amssymb}
\usepackage{latexsym}
\usepackage{amsthm}
\theoremstyle{plain}
\newtheorem{theorem}{Theorem}[section]

\newtheorem{remark}{Remark}[section]
\newtheorem{lemma}{Lemma}[section]

\makeatletter
\@addtoreset{equation}{section}

\usepackage{color}

\usepackage{graphicx,color}

\title[New examples of constant mean curvature hypersurfaces ]
{New examples of constant mean curvature hypersurfaces in the  sphere}
\author{Chuqi Huang and  Guoxin Wei}
\address{Chuqi Huang \\  School of Mathematical Sciences, South China Normal University,
510631, Guangzhou,  China, 2020021849@m.scnu.edu.cn}
\address{Guoxin Wei \\  School of Mathematical Sciences, South China Normal University,
510631, Guangzhou,  China, weiguoxin@tsinghua.org.cn}

\begin{document}

\maketitle

\begin{abstract}
In this paper, firstly, we show the existence of a compact embedded constant mean curvature (CMC) hypersurface $\Sigma_1$ in $\mathbb{S}^{2n}$ of the type $S^{n-1} \times S^{n-1} \times S^{1}$. Moreover, the hypersurface $\Sigma_1$ exhibits $O(n)\times O(n)$ symmetry. Secondly, we show that there exists a  compact embedded CMC-hypersurface $\Sigma_2 \subset \mathbb{S}^{3n-1}$ of the type $S^{n-1} \times S^{n-1} \times S^{n-1} \times S^{1}$. These results generalize the results of Carlotto and Schulz \cite{CA}.
\end{abstract}

\section{Introduction}

\noindent CMC-hypersurfaces are submanifolds with constant mean curvature H and codimension 1, which are natural extension of the minimal hypersurfaces($H=0$).
The study of CMC-hypersurfaces in the space forms is important. There are  a few examples of complete embedded CMC-hypersurfaces in the sphere.
The simplest examples of minimal hypersurfaces in $\mathbb{S}^{n+1}(1)$ are the totally geodesic n-dimensional spheres  $\mathbb{S}^{n}(1)$ and the Clifford torus $S^{k}(\sqrt{\frac{k}{n}}) \times S^{n-k}(\sqrt{\frac{n-k}{n}})$, $1\leq k \leq n-1$.
In \cite{H2}, Hsiang and Hsiang claimed that the methods developed there would allow one to prove the existence of minimal embeddings of the type $S^{1} \times S^{n-1} \times S^{n-1} \times S^{n-1}$ in the round sphere $\mathbb{S}^{3n-1}$ for any $n \geq 2$.
In 2021, Carlotto and Schulz \cite{CA} showed the existence of minimal embeddings of the type  $S^{n-1} \times S^{n-1} \times S^{1}$ in the round sphere $\mathbb{S}^{2n}$ by the method similar to that of \cite{HS} and the shooting method similar to that of \cite{AN} .
The aim of this paper is to construct some new compact, embedded CMC-hypersurfaces in $\mathbb{S}^{2n}$ and in $\mathbb{S}^{3n-1}$.  In fact, we prove the following:

\begin{theorem}\label{theorem 1.1}
For $n \geq 2$, there exists a compact embedded CMC-hypersurface $\Sigma_1 \subset \mathbb{S}^{2n}$ of the type $S^{n-1} \times S^{n-1} \times S^{1}$ and exhibits $O(n)\times O(n)$ rotational symmetry.
\end{theorem}

\begin{remark}
When $H=0$, our Theorem \ref{theorem 1.1} reduces to the reuslt of Carlotto and Schulz in \cite{CA}.
\end{remark}

\noindent To prove Theorem \ref{theorem 1.1},   we use a method similar to that of \cite{CA}.  Firstly, we consider a suitable group action on $\mathbb{S}^{2n} \subset \mathbb{R}^{2n+1}$ and study the equations obtained in the quotient space ${\mathbb{S}^{2n}/G}$; specifically, $G = O(n)\times O(n)$ acts on $\mathbb{S}^{2n}$ via the representation $\rho_n \oplus \rho_n \oplus \mathbf{1}$ as Hsiang \cite{HS}. The group action determines a relationship between the CMC-hypersurface we are interested in, and a  ''generating curve''  in $\mathbb{S}^{2n}/G$ that satisfies a system of ODEs.
Then we construct a smooth, compact generating curve that can generate a CMC-hypersurface of the type $S^{n-1} \times S^{n-1} \times S^{1}$ in the unit sphere $\mathbb{S}^{2n}$, by using a suitable shooting method  \cite{AN}.

\noindent For the case in $\mathbb{S}^{3n-1}$, we also construct the CMC-hypersurfaces $S^{n-1} \times S^{n-1} \times S^{n-1} \times S^{1}$ in $\mathbb{S}^{3n-1}$ for any $n \geq 2$.

\begin{theorem}\label{theorem 1.2}
For $n \geq 2$,  there exists a compact embedded CMC-hypersurface $\Sigma_2 \subset \mathbb{S}^{3n-1}$ of the type $S^{n-1} \times S^{n-1} \times S^{n-1} \times S^{1}$ and exhibits $O(n)\times O(n)\times O(n)$ rotational symmetry.
\end{theorem}

\section{Preliminaries}

\noindent Let $n \ge 2$, the group $G = O(n)\times O(n)$  acts on the unit sphere $\mathbb{S}^{2n}(1) = \{(\vec X, \vec Y, Z) : \vec X,\, \vec Y \in \mathbb{S}^{n-1},\, Z \in [-1, 1]\} $  in the usual way.
The quotient $\mathbb{S}^{2n}/G  = \{ (x, y, z) \in S^{2}(1) :  x = |\vec X| =\sin(r)\cos(\theta) \geq 0,\, y = |\vec Y| = \sin(r)\sin(\theta) \geq 0,\, z = Z = \cos(r)\}$ is equipped with spherical coordinates $(r, \theta) \in [0, \pi] \times [0, \frac{\pi}{2}]$ and respects to the metric
\begin{equation}
g = dr^{2} + (\sin r )^{2} d\theta^{2}.
\end{equation}

\noindent Let  $C: [0, L] \longrightarrow {(\mathbb{S}^{2n}/G)}^{\circ},\,s \mapsto \gamma (s)$ be a generating curve with arc length as parameter. And we take $\varphi (t_1,\dots\dots,t_{n-1}) = (\varphi_1,\dots \dots, \varphi_n)$ and  $\psi (u_1,\dots\dots,u_{n-1}) = (\psi_1,\dots \dots, \psi_n)$ as two orthogonal parametrizations of the unit sphere $\mathbb{S}_{n-1}(1)$. It follows that
\begin{equation}
\begin{aligned}
f(t_1,\dots\dots,t_{n-1}, u_1,\dots\dots,u_{n-1}, s) =& (x(s)\varphi_1,\dots \dots, x(s)\varphi_n, y(s)\psi_1,\dots \dots, y(s)\psi_n,z(s)),\\
\varphi_{\hat{i}}=\varphi_{\hat{i}}(t_1,\dots\dots,t_{n-1}),&\quad \varphi_1^2+ \cdots\cdots +\varphi_n^2=1,\\
\psi_{\hat{j}}=\psi_{\hat{j}}(u_1,\dots\dots,u_{n-1}),&\quad \psi_1^2+ \cdots\cdots + \psi_n^2=1,  \quad 1 \leq \hat{i},\hat{j},\hat{k}, \dots \leq n \label{f}
\end{aligned}
\end{equation}
is a parametrization of a hypersurface $\Sigma_1$ generated by the generating curve $\gamma (s)=(x(s),y(s),z(s))$, where the curve can be rewritten as $(r(s), \theta(s))$.

\noindent Since the curve belongs to ${(\mathbb{S}^{2n}/G)}^{\circ}$ and the parameter s can be chosen as its arc length, we have
\begin{displaymath}
x^2(s)+y^2(s)+z^2(s)=1, \quad \dot x^2(s)+\dot y^2(s)+\dot z^2(s)=1.
\end{displaymath}
Moreover, we can write
\begin{eqnarray}
\frac{dr}{ds} = \cos (\alpha),\label{r}\\
\nonumber\\
\frac{d\theta}{ds} = \frac{\sin (\alpha)}{\sin (r)}, \label{theta}
\end{eqnarray}
where $\alpha$ is the angle between the unit tangent vector $d\gamma/ds$ and $\partial/\partial r$.

\noindent Since, from (\ref{f}),
\begin{displaymath}
\begin{aligned}
\frac{\partial f}{\partial s}&=(\dot x(s)\varphi_1,\dots \dots, \dot x(s)\varphi_n, \dot y(s)\psi_1,\dots \dots, \dot y(s)\psi_n, \dot z(s)),\\
\frac{\partial f}{\partial t_i}&=(x(s)\frac{\partial \varphi_1}{\partial t_i},\dots \dots, x(s)\frac{\partial \varphi_n}{\partial t_i},\overbrace{0,\dots \dots ,0}^n,0),\\
\frac{\partial f}{\partial u_j}&=(\overbrace{0,\dots \dots ,0}^n, y(s)\frac{\partial \psi_1}{\partial u_j},\dots \dots, y(s)\frac{\partial \psi_n}{\partial u_j}, 0), \quad 1 \leq i, j, k, \dots \leq n-1
\end{aligned}
\end{displaymath}
and  $\varphi, \psi$ are orthogonal parametrizations,  we obtain
\begin{displaymath}
\begin{aligned}
\langle \frac{\partial f}{\partial s}, \frac{\partial f}{\partial s}\rangle &=& 1,
\qquad \langle \frac{\partial f}{\partial s},\frac{\partial f}{\partial t_i}\rangle &=& 0,
\qquad \langle \frac{\partial f}{\partial t_i},\frac{\partial f}{\partial t_j}\rangle &=& x^2(s) \sum_{\hat{k}=1}^n \left(\frac{\partial \varphi_{\hat k}}{\partial t_i}\frac{\partial \varphi_{\hat k}}{\partial t_j}\right), \\
\langle \frac{\partial f}{\partial t_i},\frac{\partial f}{\partial u_j} \rangle &=& 0,
\qquad \langle \frac{\partial f}{\partial s},\frac{\partial f}{\partial u_j} \rangle &=& 0,
\qquad \langle \frac{\partial f}{\partial u_i},\frac{\partial f}{\partial u_j} \rangle &=& y^2(s) \sum_{\hat{k}=1}^n \left(\frac{\partial \psi_{\hat k}}{\partial u_i}\frac{\partial \psi_{\hat k}}{\partial u_j}\right).
\end{aligned}
\end{displaymath}
Thus
\begin{displaymath}
\begin{aligned}
\uppercase\expandafter{\romannumeral1} &=x^2(s)\sum_{i,j}\left(\sum_{\hat{k}=1}^n \left(\frac{\partial \varphi_{\hat k}}{\partial t_i}\frac{\partial \varphi_{\hat k}}{\partial t_j}\right)dt_i dt_j \right) +  y^2(s)\sum_{i,j}\left( \sum_{\hat{k}=1}^n \left(\frac{\partial \psi_{\hat k}}{\partial u_i}\frac{\partial \psi_{\hat k}}{\partial u_j}\right)du_i du_j\right)+ ds^2 \\
&=x^2(s)d\varphi^2 + y^2(s)d\psi^2 + ds^2.
\end{aligned}
\end{displaymath}

\noindent A unit normal field N is given by
\begin{displaymath}
N = ((y\dot{z}-\dot{y}z)\varphi, (z\dot{x}-\dot{z}x)\psi, (x\dot{y}-\dot{x}y))
\end{displaymath}
and
\begin{displaymath}
\begin{aligned}
\frac{\partial^2 f}{\partial s^2}&=(\ddot x \varphi, \ddot y \psi, \ddot z),
\qquad \quad \frac{\partial^2 f}{\partial s\partial t_i}&=(\dot x \frac{\partial \varphi}{\partial t_i},\overbrace{0,\dots \dots ,0}^n, 0),
\qquad \frac{\partial^2 f}{\partial t_i\partial t_j}&=(x \frac{\partial^2 \varphi}{\partial t_i\partial t_j},\overbrace{0,\dots \dots ,0}^n, 0),\\
\frac{\partial^2 f}{\partial t_i\partial u_j}&=(\overbrace{0,\dots \dots ,0}^{2n+1}),
\qquad \frac{\partial^2 f}{\partial s\partial u_i}&=(\overbrace{0,\dots \dots ,0}^n, \dot y \frac{\partial \psi}{\partial u_i},0),
\qquad \frac{\partial^2 f}{\partial u_i\partial u_j}&=(\overbrace{0,\dots \dots ,0}^n, y \frac{\partial^2 \psi}{\partial u_i\partial u_j},0).
\end{aligned}
\end{displaymath}
Thus we have
\begin{displaymath}
\uppercase\expandafter{\romannumeral2} = \kappa_1\left(x^2(s)d\varphi^2\right) +\kappa_n\left(y^2(s)d\psi^2\right)+\kappa_{2n-1}ds^2,
\end{displaymath}
where the principal curvatures
\begin{equation}
\begin{aligned}
\kappa_1 = \dots =\kappa_{n-1}  =-\frac{y\dot{z}-\dot{y}z}{x}=\frac{\cos r \cos \theta \sin \alpha + \sin \theta \cos \alpha }{\sin r \cos \theta}, \label{k1}
\end{aligned}
\end{equation}
\begin{equation}
\begin{aligned}
\kappa_n = \dots = \kappa_{2n-2} =-\frac{z\dot{x}-\dot{z}x}{y}=\frac{\cos r\sin \theta\sin \alpha - \cos \theta \cos \alpha}{\sin r \sin \theta}, \label{k2}
\end{aligned}
\end{equation}
\begin{equation}
\begin{aligned}
\kappa_{2n-1} = (y\dot{z}-\dot{y}z)\ddot{x}+(z\dot{x}-\dot{z}x)\ddot{y}+(x\dot{y}-\dot{x}y)\ddot{z}=\frac{d\alpha} {ds}+ \cot r \sin \alpha. \label{k3}
\end{aligned}
\end{equation}

\noindent From (\ref{k1}), (\ref{k2}), (\ref{k3}), the CMC-hypersurface equation with $H = \sum_{i=1}^{2n-1}\kappa_{i} \equiv \lambda$ reduces to
\begin{equation}
\frac{d\alpha} {ds} = (2n-2)\frac{\cot (2\theta)}{\sin r} \cos \alpha - (2n-1)\cot r \sin \alpha +  \lambda. \label{alpha}
\end{equation}

\noindent We  rewrite the differential equation above as an $3 \times 3$ ODE system, i.e. to
have it in the form $dU/ds = f(U)$ where $U$ varies in a suitable subset of $\mathbb{R}^3$. Indeed, we can take $U = (r, \theta, \alpha)$ and thus consider the system of equations
\begin{displaymath}
\left\{ \begin{array}{ll}
\frac{dr}{ds} = \cos (\alpha),  & \textrm{(\ref{r})}\\
\, \\
\frac{d\theta}{ds} = {\sin (\alpha)}/{\sin (r)},  & \textrm{(\ref{theta})}\\
\, \\
\frac{d\alpha} {ds} = (2n-2)\frac{\cot (2\theta)}{\sin (r)} \cos(\alpha) - (2n-1)\cot(r)\sin(\alpha) +  \lambda.  & \textrm{(\ref{alpha})}
\end{array} \right.
\end{displaymath}

\newcommand{\arccot}{\mathrm{arccot}\,}

\noindent We consider two special cases of the above ODE system:

\noindent{\bf{Case (a):}} Assuming r = const, that is $\frac{dr}{ds} = \cos (\alpha) \equiv 0 $, then we get special solutions:
\begin{displaymath}
\begin{aligned}
&\text{when} \quad \alpha \equiv \frac{\pi}{2}, \qquad \frac{d\alpha} {ds} \equiv 0 \Rightarrow  r \equiv \arccot \left(\frac{\lambda}{2n-1} \right);\\
&\text{when} \quad \alpha \equiv -\frac{\pi}{2}, \qquad \frac{d\alpha} {ds} \equiv 0 \Rightarrow  r \equiv \arccot \left(\frac{-\lambda}{2n-1} \right).
\end{aligned}
\end{displaymath}
The solutions correspond to the $(2n-1)$-dimensional spheres in $\mathbb{S}^{2n}$ respectively.

\noindent{\bf{Case (b):}} Assuming  $\theta$ = const, that is $\frac{d\theta}{ds} = \frac{\sin (\alpha)}{\sin (r)} \equiv 0 $. When $\lambda \neq 0$, there is no solution that satisfies this assumption, which is different from the case of minimal hypersurfaces.
\begin{displaymath}
\text{When} \quad \alpha \equiv 0, \qquad \frac{d\alpha} {ds} \equiv 0 (\lambda \neq 0) \Rightarrow  r = const \Rightarrow \frac{dr}{ds} \equiv 0,\qquad \text{contradict with } \frac{dr}{ds} = \cos (\alpha) = 1.
\end{displaymath}

\noindent In Section 3, we will use a suitable shooting method to construct solutions of the system (\ref{r}), (\ref{theta}), (\ref{alpha}) such that the projections of the trajectories on the  quotient $\mathbb{S}^{2n}/G$ are simple and periodic. Actually, these curves generate the embedded CMC-hypersurfaces claimed in the statement of Theorem \ref{theorem 1.1}.

\section{Proof of Theorem Theorem \ref{theorem 1.1}}
\noindent {\bf{Symmetry of the solutions}}. \quad
Let $(r(s), \theta(s), \alpha(s))$ be a solution of the above system corresponding to  a curve $C: [0, L] \longrightarrow {(\mathbb{S}^{2n}/G)}^{\circ}$. It is easy to verify that the solution has symmetry about $\{\theta = \frac{\pi}{4}\}$ i.e. the symmetry of its corresponding curve about $\{\theta = \frac{\pi}{4}\}$ still corresponds to a solution of the above system: $(\tilde{r}(s), \tilde{\theta}(s), \tilde{\alpha}(s)) = (r(L - s), \frac{\pi}{2} - \theta(L - s), \pi - \alpha(L - s))$. Meanwhile, the solution also has symmetry about $\{r = \frac{\pi}{2}\}$, i.e. $(\bar{r}(s), \bar{\theta}(s), \bar{\alpha}(s)) = (\pi - r(L - s), \theta(L-s), - \alpha(L - s))$  is also a solution of the above system.

\noindent Assume $H > 0$, that is $\lambda > 0$ ( when $H < 0$, we consider the reverse curve that corresponding to the solution $(\hat{r}(s), \hat{\theta}(s), \hat{\alpha}(s)) = (r(L - s), \theta(L - s), \pi +  \alpha(L - s))$. A direct calculation shows that $\hat{H} = - H$).

\begin{lemma} Given any $r_0 \in (0, \frac{\pi}{2})$, there exists $s_\ast > 0$ such that the system (\ref{r}), (\ref{theta}), (\ref{alpha}) has a unique solution
\begin{displaymath}
(r, \theta, \alpha) : [0, s_\ast] \to B := (0, \frac{\pi}{2}] \times (0, \frac{\pi}{4}] \times [-\frac{\pi}{2}, 0]
\end{displaymath}
with initial data $(r, \theta, \alpha)(0) = (r_0, \frac{\pi}{4}, - \frac{\pi}{2})$ on which the solution depends continuously such that $r(s_\ast) = \frac{\pi}{2}$ or $\alpha(s_\ast) = 0$. Moreover, $dr/ds \geq 0$, $d\theta/ds \leq 0$, $d\alpha/ds > 0$ and $\theta(s) >0 $ for all $s \in [0, s_\ast]$.
\end{lemma}

\noindent {\bf Proof.}  As long as $(r, \theta, \alpha) \in B:= (0, \frac{\pi}{2}] \times (0, \frac{\pi}{4}] \times [-\frac{\pi}{2}, 0]$ and $\lambda > 0$, equations (\ref{r}), (\ref{theta}), (\ref{alpha}) imply $dr/ds \geq 0, d\theta/ds \leq 0$ and  $d\alpha/ds > 0$. Therefore, solutions with initial data $(r, \theta, \alpha)(0) = (r_0, \frac{\pi}{4}, - \frac{\pi}{2})$ stay in $B$ at least for a short time and can leave $B$ only via $r = \frac{\pi}{2}$, $\theta = 0$ or $\alpha = 0$.

\noindent To show that $\theta(s_\ast)=0$ can not occur, let's assume there exists a trajectory leaves $B$ via $\theta = 0$, and let $s_1 \in (0, s_\ast)$ be fixed such that $\theta(s_1) > \frac{\pi}{8}$. Equations (\ref{theta}) and (\ref{alpha}) imply
\begin{displaymath}
\begin{aligned}
\frac{d\alpha} {ds}  & = [(2n-2)\cot(2\theta)\cot(\alpha) - (2n-1)\cos(r)] \frac{d\theta}{ds} + \lambda \\
& \geq (2n-2)\cot(2\theta)\cot(\alpha_1)\frac{d\theta}{ds}
\end{aligned}
\end{displaymath}
for all $s \in [s_1, s_\ast)$. Integrating this differential inequality over an interval $[s_1, s)$ gives
\begin{displaymath}
\frac{\pi}{2} > \alpha(s)\arrowvert
_{s_1}^{s} \geq (n-1)\cot(\alpha_1) \ln (\sin(2\theta(s)))\arrowvert_{s_1}^{s}.
\end{displaymath}
When $s \to s_{\ast}$, the above inequality implies that $\frac{\pi}{2} > \infty$, which gives a contradiction to the fact. Thus the trajectory in $B$ must stay at positive distance away from $\{\theta = 0\}$.

\noindent In particular, since the vector field (cf. right-hand side of (\ref{r}), (\ref{theta}), (\ref{alpha})) is smooth and uniformly bounded in any closed subset of $B$, from the continuous dependence of the solutions on initial data, we have, for a given $r_0 \in (0, \pi/2)$ and any $\varepsilon > 0$, there exists a sufficiently small $\eta > 0$ such that the image through the flow map of $[r_0 - \eta, r_0 + \eta]$ will stay at positive distance from the plane$\{\theta = 0\}$, and  these trajectories will be $\varepsilon$-close at any time (in $B$).
Moreover, from $\frac{dr}{ds}(s_\ast)>0$ and $\frac{d\alpha}{ds}(s_\ast)>0$, we can see that these trajectories with initial data in $[r_0 - \eta, r_0 + \eta]$ will leave $B$ at times less than $\varepsilon$ apart.
\begin{flushright}
Q. E. D.
\end{flushright}

\noindent To simplify the calculation, we consider the following two cases.

\noindent {\bf Case 1}: $\lambda \in (\frac{4}{\pi}, \infty).$
\begin{lemma}  For $n \geq 2$ and $\lambda  \in (\frac{4}{\pi}, \infty)$, if the initial data $r_0 \in (0, \frac{\pi}{2})$ is sufficiently small, then the solution constructed in Lemma 3.1 satisfies $\alpha(s_\ast) = 0.$
\end{lemma}
\noindent {\bf Proof.}  From equations (\ref{r}) and (\ref{alpha}), we obtain
\begin{equation}
\begin{aligned}
\frac{d\alpha} {ds}  & = [(2n-2)\frac{\cot(2\theta)}{\sin(r)} - (2n-1)\cot(r)\tan(\alpha)+\frac{\lambda}{\cos(\alpha)}]\frac{dr}{ds}\\
& \geq \frac{\lambda}{\cos(\alpha)}\frac{dr}{ds} \label{9}
\end{aligned}
\end{equation}
for all $s \in (0, s_\ast]$. The inequality (\ref{9}) implies
\begin{equation}
\cos(\alpha) \frac{d\alpha}{ds} \geq \lambda \frac{dr}{ds}.\label{10}
\end{equation}
Since  $\lambda > \frac{4}{\pi}$, from (\ref{10}), we then have
\begin{displaymath}
\begin{aligned}
&1 \geq \sin(\alpha(s))\arrowvert_{0}^{s_\ast} \geq \lambda r(s)\arrowvert
_{0}^{s_\ast} \\
& \Rightarrow \frac{\pi}{4} > \frac{1}{\lambda} \geq r(s_\ast) - r_0 \\
& \Rightarrow \frac{\pi}{4} + r_0 > r(s_\ast).
\end{aligned}
\end{displaymath}

\noindent If the initial value $r_0 < \frac{\pi}{4}$, then $r(s_\ast) < \frac{\pi}{2}$, and we obtain $\alpha(s_\ast) = 0$ from Lemma 3.1.
\begin{flushright}
Q. E. D.
\end{flushright}

\begin{lemma}  For $n \geq 2$ and $\lambda  \in (\frac{4}{\pi}, \infty)$, let $(r, \theta, \alpha) : [0, s_\ast] \to B:= (0, \frac{\pi}{2}] \times (0, \frac{\pi}{4}] \times [-\frac{\pi}{2}, 0]$ be the solution constructed in Lemma 3.1. There exists  a constant $r_\lambda \in (0, \frac{\pi}{2})$ depending only on $\lambda$  such that if the initial data $r_0 \in (r_\lambda, \frac{\pi}{2})$, then the soluiton constructed in Lemma 3.1 satisfies
\begin{displaymath}
0 \leq \cot(2\theta(s)) \leq \cot(\frac{\pi}{2} - \frac{2}{\lambda \sin(r_0)}), \qquad \forall s \in [0, s_\ast].
\end{displaymath}
\end{lemma}
\noindent {\bf Proof.} Equations (\ref{theta}) and (\ref{alpha})  show that
\begin{equation}
\begin{aligned}
\frac{d\alpha} {ds}  & = (2n-2)\frac{\cot (2\theta)}{\sin (r)} \cos(\alpha) - (2n-1)\cot(r)\sin(\alpha) +  \lambda \\
& \geq \frac{\lambda \sin(r)}{\sin(\alpha)}\frac{d\theta}{ds} \label{11}
\end{aligned}
\end{equation}
for all $s \in [0, s_\ast)$. The inequality (\ref{11}) implies
\begin{equation}
-\sin(\alpha) \frac{d\alpha}{ds} \geq - \lambda \sin(r_0) \frac{d\theta}{ds}.\label{12}
\end{equation}
Integrating (\ref{12}), we get
\begin{displaymath}
\begin{aligned}
&1 \geq \cos(\alpha(s))\arrowvert_{0}^{s_\ast} \geq -\lambda \sin(r_0)\theta(s) \arrowvert_{0}^{s_\ast} \\
& \Rightarrow \theta(s_\ast) \geq \frac{\pi}{4} - \frac{1}{\lambda \sin(r_0)} .
\end{aligned}
\end{displaymath}
Since  $\lambda > \frac{4}{\pi}$, there exists $r_\lambda$ sufficiently close to $\frac{\pi}{2}$ such that $\frac{1}{\lambda \sin(r_\lambda)} < \frac{\pi}{4}.$
Therefore, given any $r_0 \in (r_\lambda, \frac{\pi}{2})$, we have $\theta(s) \geq \theta(s_\ast) \geq \frac{\pi}{4} - \frac{1}{\lambda \sin(r_0)} > 0$ for all $s \in [0, s_\ast], $ that is  $0 \leq \cot(2\theta(s)) \leq \cot(\frac{\pi}{2} - \frac{2}{\lambda \sin(r_0)}).$

\begin{flushright}
Q. E. D.
\end{flushright}

\begin{lemma}  For $n \geq 2$ and $\lambda  \in (\frac{4}{\pi}, \infty)$, if the initial data $r_0 \in (0, \frac{\pi}{2})$ is sufficiently close to $\frac{\pi}{2}$, then the solution constructed in Lemma 3.1 satisfies $r(s_\ast) = \frac{\pi}{2}.$
\end{lemma}
\noindent {\bf Proof.}  If $\alpha(s_\ast) \neq 0$ for all $r_0 \in (r_\lambda, \frac{\pi}{2})$, then we must have $r_\ast = \frac{\pi}{2}$ for all $r_0 \in (r_\lambda, \frac{\pi}{2})$ and the Lemma is proved.

\noindent We suppose there is a convergent sequence $\{r_n\}$ which converges to $\frac{\pi}{2}$ such that the solutions constructed in Lemma 3.1  with initial data $r(0)=r_n$ satisfy $\alpha(s_\ast) = 0$. Then we can choose $s_1 \in (0, s_\ast)$ such that $\alpha_1 = \alpha(s_1) = -\frac{3\pi}{8}$. When the initial data $r_n > r_\lambda$, we have
\begin{equation}
\begin{aligned}
\frac{d\alpha} {ds} & = [(2n-2)\frac{\cot(2\theta)}{\sin(r)} - (2n-1)\cot(r)\tan(\alpha)+\frac{\lambda}{\cos(\alpha)}]\frac{dr}{ds}\\
& \leq [(2n-2)\frac{\cot(\frac{\pi}{2} - \frac{2}{\lambda \sin(r_n)})}{\sin(r_n)} - (2n-1)\cot(r_n)\tan(\alpha_1)+\frac{\lambda}{\cos(\alpha_1)}]\frac{dr}{ds} \label{13}
\end{aligned}
\end{equation}
for all $s \in [s_1, s_\ast]$. We obtain the following inequality by integrating  (\ref{13}),
\begin{equation}
\frac{3\pi}{8} = \alpha(s) \arrowvert_{s_1}^{s_\ast} \leq  [(2n-2)\frac{\tan(\frac{2}{\lambda \sin(r_n)})}{\sin(r_n)} - (2n-1)\cot(r_n)\tan(\alpha_1)+ \frac{\lambda}{\cos(\alpha_1)}] r(s) \arrowvert_{s_1}^{s_\ast}. \label{14}
\end{equation}
When $r_n \to \frac{\pi}{2}$, inequality (\ref{14}) implies
\begin{displaymath}
\frac{3\pi}{8} = \alpha(s) \arrowvert_{s_1}^{s_\ast} \leq 0,
\end{displaymath}
which gives a contradiction to the fact. According to Lemma 3.1, we have $r(s_\ast) = \frac{\pi}{2}$ when $r_0$ is sufficiently close to $\frac{\pi}{2}$.

\begin{flushright}
Q. E. D.
\end{flushright}

\noindent {\bf Case 2}: $\lambda \in (0, \frac{4}{\pi}].$
\begin{lemma}  For $n \geq 2$ and $\lambda  \in (0, \frac{4}{\pi}]$, given $r_0 \in (0, \frac{\pi}{4})$, let $(r, \theta, \alpha) : [0, s_\ast] \to B:= (0, \frac{\pi}{2}] \times (0, \frac{\pi}{4}] \times [-\frac{\pi}{2}, 0]$ be the solution constructed in Lemma 3.1 and let $s_1 \in [0, s_\ast]$ be arbitrary. If $r(s_1) \geq 2r_0$,  then
\begin{displaymath}
\theta (s_1) < \frac{\pi}{4} - \frac{1}{6n}.
\end{displaymath}
\end{lemma}
\noindent {\bf Proof.}  Let $\delta = \frac{1}{6n}$ and assume there exists $s_1 \in (0, s_\ast]$ with $r(s_1) \geq 2 r_0$ (Otherwise, the lemma holds directly). We suppose $\theta(s_1) \geq \frac{\pi}{4} - \frac{1}{6n}$. By monotonicity of $\theta$,  we then have  $\theta(s) \geq \frac{\pi}{4} - \delta, \forall s \in [0, s_1]$. Hence,
\begin{displaymath}
\theta(s) \geq \frac{\pi}{4} - \delta,\quad 0 \leq \cot(2\theta) \leq  \tan(2\delta) := b, \quad \cot(\alpha) \leq 0,\quad \sin(\alpha)<0, \quad 0 \leq \cos(r) \leq 1, \quad \frac{d\theta}{ds} \leq 0
\end{displaymath}
for all $s \in [0, s_1)$. Since $\lambda \leq \frac{4}{\pi} <2$, equations (\ref{theta}) and (\ref{alpha}) imply
\begin{equation}
\begin{aligned}
\frac{d\alpha} {ds}  & = [(2n-2)\cot(2\theta)\cot(\alpha)- (2n-1)\cos(r) + \frac{\lambda\sin(r)}{\sin(\alpha)}]\frac{d\theta}{ds}\\
\,
& \leq (2n-1)(b\cot(\alpha) - 1 + \frac{1}{\sin(\alpha)})\frac{d\theta}{ds} \label{15}
\end{aligned}
\end{equation}
for all $s \in [0, s_1)$. The inequality (\ref{15}) gives us
\begin{equation}
\frac{\sin(\alpha)}{b\cos(\alpha) - \sin(\alpha)+ 1} \frac{d\alpha} {ds}\geq (2n-1)\frac{d\theta}{ds}. \label{16}
\end{equation}
Integrating the inequality (\ref{16}) over $[0, s_1)$ yields
\begin{displaymath}
\begin{aligned}
 - \frac{1}{b^2+1}[&b\ln(b\cos(\alpha) - \sin(\alpha)+ 1) + \alpha + \frac{1}{b} \ln(\frac{\tan(\alpha/2) - 1}{\tan(\alpha/2) - \frac{1+b}{1-b}})] \arrowvert_{-\frac{\pi}{2}}^{\alpha(s_1)} \geq (2n-1)\theta \arrowvert_{\frac{\pi}{4}}^{\theta(s_1)} \geq (2n-1)(-\delta) \\
& \Rightarrow (b^2+1)(2n-1)\delta \geq [b\ln(b\cos(\alpha) - \sin(\alpha)+ 1) + \alpha + \frac{1}{b} \ln(\frac{\tan(\alpha/2) - 1}{\tan(\alpha/2) - \frac{1+b}{1-b}})] \arrowvert_{-\frac{\pi}{2}}^{\alpha(s_1)}.
\end{aligned}
\end{displaymath}
Defining the functions $G(\alpha):=b\ln(b\cos(\alpha) - \sin(\alpha)+ 1) + \alpha$ and $F(y):=\frac{1}{b} \ln(\frac{y+1}{y+\frac{1+b}{1-b}})$, where $y:=-\tan(\alpha/2)$. Then we have
\begin{displaymath}
\begin{aligned}
&G(\alpha)=b\ln(b\cos(\alpha) - \sin(\alpha)+ 1) + \alpha \geq  b\ln(b+1)+\alpha, \quad \forall \alpha \in [-\frac{\pi}{2}, 0]\\
&\Rightarrow G(\alpha)-G(-\frac{\pi}{2}) \geq b\ln(\frac{b+1}{2})+\alpha+\frac{\pi}{2}.
\end{aligned}
\end{displaymath}
The map $f: [0,1] \to \mathbb{R}$ given by
\begin{displaymath}
f(y)= \ln(\frac{y+1}{(1-b)y+1+b})
\end{displaymath}
is increasing with minimum value $f(0)=-\ln(1+b)$, hence we obtain
\begin{displaymath}
\begin{aligned}
&F(y)-F(-\tan(-\frac{\pi}{4}))=\frac{1}{b} (\ln(\frac{y+1}{y+\frac{1+b}{1-b}})-\ln(1-b)) = \frac{1}{b}\ln(\frac{y+1}{(1-b)y+1+b})\\
&\Rightarrow F(-\tan(\alpha/2))-F(-\tan(-\frac{\pi}{4})) \geq -\frac{1}{b}\ln(1+b)
\end{aligned}
\end{displaymath}
for all $\alpha \in [-\frac{\pi}{2}, 0]$.

Now, we have
\begin{displaymath}
\begin{aligned}
(b^2+1)(2n-1)\delta &\geq G(\alpha(s_1))-G(-\frac{\pi}{2})+F(-\tan(\alpha(s_1)/2))-F(-\tan(-\frac{\pi}{4})) \\
&\geq b\ln(\frac{b+1}{2})+\alpha(s_1)+\frac{\pi}{2}-\frac{1}{b}\ln(1+b),
\end{aligned}
\end{displaymath}
which can be rewritten as
\begin{displaymath}
\alpha(s_1) \leq (b^2+1)(2n-1)\delta-b\ln(\frac{b+1}{2})-\frac{\pi}{2}+\frac{1}{b}\ln(1+b).
\end{displaymath}
Since $\delta = \frac{1}{6n}$, we obtain
\begin{displaymath}
\begin{aligned}
\alpha(s_1) &\leq (b^2+1)/3 - b\ln(\frac{b+1}{2})-\frac{\pi}{2}+\frac{1}{b}\ln(1+b)\\
&\leq - 0.15
\end{aligned}
\end{displaymath}
for all $b \in (0, 1)$. By monotonicity of $\alpha$, we have $\alpha(s) < -0.15$ for all $s \in [0, s_1]$. From (\ref{r}) and (\ref{theta}), we then get
\begin{equation}
\frac{d\theta}{ds}=\frac{\tan(\alpha)}{\sin(r)} \frac{dr}{ds} \leq \frac{\tan(-0.15)}{\sin(r)})\frac{dr}{ds}. \label{17}
\end{equation}
Integrating the inequality (\ref{17}) over $[0, s_1)$ implies
\begin{equation}
\begin{aligned}
\theta \arrowvert_{\frac{\pi}{4}}^{\theta(s_1)} & \leq \tan(-0.15) \ln(\tan(r/2))\arrowvert_{r_0}^{r(s_1)}\\
& \leq -\tan(0.15)\ln(2)\\
& < - \frac{1}{6n},
\end{aligned}
\end{equation}
which gives a contradiction to our initial assumption. Thus we prove the claim.
\begin{flushright}
Q. E. D.
\end{flushright}

\begin{lemma}  For $n \geq 2$ and $\lambda  \in (0, \frac{4}{\pi}]$, there exists a constant $c_n>0$ depending only on n such that if the initial data $r_0 \in (0, \frac{\pi}{2})$ is sufficiently small, then the solution constructed in Lemma 3.1 satisfies $\alpha(s_\ast) = 0$ and $r(s_\ast) \leq c_n r_0.$
\end{lemma}
\noindent {\bf Proof.} For a given $r_0 \in (0, \frac{\pi}{4})$, we may assume $r(s_\ast) > 2r_0$, otherwise the claim follows directly. Let $s_1 \in (0, s_\ast)$ be fixed such that $r(s_1) = 2r_0$. According to Lemma 3.5 and the monotonicity of $\theta$, we have
\begin{displaymath}
\theta (s) < \frac{\pi}{4} - \frac{1}{6n}
\end{displaymath}
for all $s\in[s_1, s_\ast]$, which implies $\cot(2\theta) > \tan(1/(3n))$ for all $s\in[s_1, s_\ast]$. Hence,
\begin{equation}
\begin{aligned}
\frac{d\alpha} {ds} & = [(2n-2)\frac{\cot(2\theta)}{\sin(r)} - (2n-1)\cot(r)\tan(\alpha)+\frac{\lambda}{\cos(\alpha)}]\frac{dr}{ds}\\
& \geq \frac{2n-2}{r}\tan(\frac{1}{3n})\frac{dr}{ds} \label{19}
\end{aligned}
\end{equation}
for all $s\in[s_1, s_\ast]$. It follows that
\begin{displaymath}
\begin{aligned}
&\frac{\pi}{2} \geq \alpha(s) \arrowvert_{s_1}^{s_\ast} \geq (2n-2)\tan\left(\frac{1}{3n}\right)\ln((r(s))) \arrowvert_{s_1}^{s_\ast} \\
& \Rightarrow r(s_\ast) \leq 2\exp(\frac{\pi}{4n-4}\cot(\frac{1}{3n}))r_0 := c_n r_0.
\end{aligned}
\end{displaymath}
If the initial value $r_0>0$ is chosen such that  $c_n r_0<\frac{\pi}{2}$, then $r(s_\ast)<\frac{\pi}{2}$. Thus we obtain $\alpha(s_\ast)=0$ by Lemma 3.1.

\begin{flushright}
Q. E. D.
\end{flushright}

\begin{lemma} For $n \geq 2$ and $\lambda  \in (0, \frac{4}{\pi}]$, if the initial data $r_0 \in (0, \frac{\pi}{2})$ is sufficiently close to $\frac{\pi}{2}$, then the solution constructed in Lemma 3.1 satisfies $r(s_\ast) = \frac{\pi}{2}.$
\end{lemma}
\noindent {\bf Proof.}  Let $x=\tan(r)$, $y =  \cot(2\theta)$,  and $z=-\cot(\alpha)$, we then have
\begin{displaymath}
\begin{aligned}
\sin(r)&=&\frac{x}{\sqrt{x^2+1}}, \qquad \qquad \sin(2\theta) &=& \frac{1}{\sqrt{y^2+1}}, \qquad \sin(\alpha)&=&-\frac{1}{\sqrt{z^2+1}},\\
\cos(r)&=&\frac{1}{\sqrt{x^2+1}}, \qquad  \qquad\cos(2\theta) &=& \frac{y}{\sqrt{y^2+1}}, \qquad \cos(\alpha)&=&\frac{z}{\sqrt{z^2+1}}.
\end{aligned}
\end{displaymath}
For  $s \in [0, s_\ast)$, the system (\ref{r}), (\ref{theta}), (\ref{alpha}) has the equivalent form
\begin{eqnarray}
\frac{dx}{ds} &=& \frac{(x^2+1)z}{\sqrt{z^2+1}},\label{20} \\
\nonumber\\
\frac{dy}{ds} &= &\frac{2(y^2+1)\sqrt{x^2+1}}{x\sqrt{z^2+1}}, \label{21} \\
\nonumber\\
\frac{dz} {ds} &= &(z^2+1)[(2n-2)\frac{yz\sqrt{x^2+1}}{x\sqrt{z^2+1}} +(2n-1) \frac{1}{x\sqrt{z^2+1}} +  \lambda]\label{22}
\end{eqnarray}
with initial values $x(0) = x_0 := \tan(r_0)>0$ and $y(0)=0=z(0)$. Since the domain of the system above is $(x,y,z) \in Q := (0,\infty) \times [0, \infty) \times [0, \infty)$, we have $\frac{dx}{ds} \ge 0$, $\frac{dy}{ds} \ge 0$ and $\frac{dz}{ds}\ge 0$ for all $s \in [0, s_\ast]$.

\noindent Equations (\ref{21}) and (\ref{22}) imply
\begin{equation}
\frac{dz} {ds} = [\frac{(n-1)yz(z^2+1)}{y^2+1}+ \frac{(2n-1)(z^2+1)}{2(y^2+1)\sqrt{x^2+1}}] \frac{dy}{ds}+ \lambda(z^2+1). \label{23}
\end{equation}
For any given $s_1 \in (0, s_\ast)$, let $x_1=x(s_1), y_1=y(s_1),z_1=z(s_1)$ and
\begin{displaymath}
\varepsilon = \frac{(2n-1)}{2\sqrt{x_1^2+1}}.
\end{displaymath}
Since $\frac{dz}{ds}(0)>\varepsilon \frac{dy}{ds}(0) + \lambda$, there exists $s_\varepsilon \in (0, s_1]$ such that $z \geq \varepsilon y+\lambda s$ for all $s \in [0,s_\varepsilon]$. In particular,
\begin{equation}
\begin{aligned}
\frac{dz} {ds} &\geq (z^2+1)[\frac{(n-1)\varepsilon y^2}{y^2+1}+ \frac{\varepsilon}{(y^2+1)}] \frac{dy}{ds}+ \lambda(z^2+1)\\
&\geq (z^2+1)[\varepsilon \frac{dy}{ds} + \lambda] \label{24}
\end{aligned}
\end{equation}
for all $s \in [0,s_\varepsilon]$. The inequality (\ref{24}) can be written as
\begin{equation}
\frac{1}{z^2+1} \frac{dz} {ds} \geq\varepsilon \frac{dy}{ds} + \lambda \label{25}
\end{equation}
for all $s \in [0,s_\varepsilon]$. Integrating (\ref{25}) over $[0, s_\varepsilon]$, we then have
\begin{equation}
\begin{aligned}
\arctan(z(s) \arrowvert_{0}^{s_\varepsilon} &\geq [\varepsilon y(s)+\lambda s ]\arrowvert_{0}^{s_\varepsilon} \\
\Rightarrow z(s_\varepsilon) > \tan(\varepsilon y(s_\varepsilon)&+\lambda s_\varepsilon) >\varepsilon y(s_\varepsilon)+\lambda s_\varepsilon. \label{26}
\end{aligned}
\end{equation}
Thus $s_\varepsilon = s_1$. Moreover, we get
\begin{equation}
y_1 \leq \frac{z_1-\lambda s_1}{\varepsilon} \leq (z_1 - \lambda s_1) \sqrt{x_1^2+1}. \label{27}
\end{equation}

\noindent Since equation (\ref{23}) shows
\begin{equation}    	
\frac{d}{ds}(\ln(\frac{z^2}{z^2+1})) = \frac{2}{z(z^2+1)}\frac{dz}{ds} = [\frac{(2n-2)y}{y^2+1}+\frac{(2n-1)}{(y^2+1)z\sqrt{x^2+1}}]\frac{dy}{ds}+\frac{2\lambda}{z} \geq (n-1)\frac{d}{ds}(\ln(y^2+1))+\frac{2\lambda}{z} \label{28}
\end{equation}
for all $s\in(0,s_\ast)$, we obtain
\begin{displaymath}
\begin{aligned}
 &\qquad \qquad \qquad \quad \ln(\frac{z^2(s)}{z^2(s)+1})\arrowvert_{s_1}^{s} &\geq& \qquad [(n-1)\ln(y^2(s)+1))+\frac{2\lambda}{z(s)}s]\arrowvert_{s_1}^{s}\\
&\Rightarrow \qquad \qquad \quad \frac{z(s)^2}{z(s)^2+1} \cdot \frac{z_1^2+1}{z_1^2} &\geq& \qquad  (n-1)\frac{y^2(s)+1}{(y_1^2+1)} \cdot \exp(\frac{2\lambda(s - s_1)}{z(s)})
\end{aligned}
\end{displaymath}
for all $s\in[s_1,s_\ast)$. Combining (\ref{27}), it follows that
\begin{equation}
\frac{(n-1)y^2(s)(z^2(s)+1) \exp(\frac{2\lambda (s - s_1)}{z(s)})}{z^2(s)} < \frac{z_1^2+1}{z_1^2}(\frac{4(x_1^2+1)z_1^2}{(2n-1)^2}+1). \label{29}
\end{equation}

\noindent From (\ref{20}), (\ref{22}) and (\ref{29}), we can see
\begin{equation}
\begin{aligned}
\frac{dz}{ds} &=\frac{(1+z^2)^{3/2}}{(x^2+1)z}[(2n-2)\frac{yz\sqrt{x^2+1}}{x\sqrt{z^2+1}}+\frac{2n-1}{x\sqrt{z^2+1}}+\lambda]\frac{dx}{ds}\\
&=[(2n-2)\frac{y\sqrt{z^2+1}}{z\sqrt{x^{-2}+1}}+\frac{(2n-1)\sqrt{z^2+1}}{x(x^{-2}+1)z^2}+\frac{1+z^2}{z^2(1+x^{-2})}\lambda]\frac{z\sqrt{z^2+1}}{x^2}\frac{dx}{ds}\\
&\leq [(2n-2)\sqrt{1+z_1^{-2}}\sqrt{\frac{4(x_1^2+1)z_1^2}{(2n-1)^2}+1}+\frac{(2n-1)\sqrt{z_1^2+1}}{x_1 z_1^2}+\lambda(1+z_1^{-2})] \frac{z\sqrt{z^2+1}}{x^2}\frac{dx}{ds} \\
&: =C_1 \frac{z\sqrt{z^2+1}}{x^2}\frac{dx}{ds}  \label{30}
\end{aligned}
\end{equation}
for all $s \in [s_1,s_\ast)$, where we have used the fact that the function $0<z \mapsto z^{-2}\sqrt{z^2+1}$ is decreasing.
Dividing  (\ref{30}) by $z\sqrt{z^2+1}$ and integrating, it is clear that
\begin{equation}
\int_{s_1}^{s} \frac{dz/ds}{z\sqrt{z^2+1}}ds \leq \frac{C_1}{x(s)}\arrowvert_{s_1}^{s} \leq \frac{C_1}{x_1} \label{31}
\end{equation}
for all $s \in [s_1,s_\ast).$

\noindent We suppose that $s \to z(s)$ is unbounded for any choice of initial value $x_0 >1$. Then, we can choose $s_1 \in (0,s_\ast)$ such that $z_1 = \frac{1}{\sqrt{x_0}}$. With this choice, the inequality (\ref{31}) implies
\begin{equation}
\int_{\frac{1}{\sqrt{x_0}}}^{1} \frac{1}{z\sqrt{z^2+1}}dz < \frac{C_1}{x_1}.\label{32}
\end{equation}
When $x_0 \to \infty$, the inequality (\ref{32}) implies $+ \infty < \frac{4n-4}{2n-1}+\lambda$. This contradiction shows that $z$ is bounded. Thus $\alpha(s_\ast)$ with the initial data $r_0$ close enough to $\frac{\pi}{2}$ is negative, which means $r(s_\ast)=\frac{\pi}{2}$.
\begin{flushright}
Q. E. D.
\end{flushright}

\begin{remark}
The inequality (\ref{26}) shows that $\varepsilon y(s) + \lambda s  <  \frac{\pi}{2}$ for all $s \in [0, s_\ast)$ by arbitrariness of $s_1$.
Thus we have $s_\ast \leq \frac{\pi}{2\lambda}$, which means that the length of the compact generating curve we constructed has a upper bound depending only on $\lambda$, i.e. $L(C) \leq \frac{2\pi}{\lambda}$.
\end{remark}

\begin{remark}
In fact, Lemma 3.7 is true for all  $\lambda>0$.
\end{remark}

\noindent {\bf Proof of Theorem \ref{theorem 1.1}}. Now we have found that for $n \geq 2$ and $\lambda  > 0$, there exist $r_0',\,r_0'' \in (0, \frac{\pi}{2})$ with $r_0' < r_0''$ such that the trajectories emanating from $(r_0', \frac{\pi}{4}, -\frac{\pi}{2})$ and $(r_0'', \frac{\pi}{4}, -\frac{\pi}{2})$ leave the domain $B := (0, \frac{\pi}{2}] \times (0, \frac{\pi}{4}] \times [-\frac{\pi}{2}, 0]$ from $\{\alpha = 0\}$ and from $\{r=\frac{\pi}{2}\}$ respectively.

\noindent Since the solutions constructed in Lemma 3.1 depend continuously on the initial data  $(r_0, \frac{\pi}{4}, - \frac{\pi}{2})$ and satisfy $dr/ds \geq 0$, $d\theta/ds \leq 0$, $d\alpha/ds > 0$  and $\theta(s) > 0$ for all $s \in [0, s_\ast]$,  there exists a $r_0 \in (r_0', r_0'')$ (by the continuity of the system of ODEs), for which, the trajectory emanating from $(r_0, \frac{\pi}{4}, -\frac{\pi}{2})$ leaves the domain $B$ from $\{\alpha= 0\} \cap \{r = \frac{\pi}{2}\}$. The projection of this trajectory on the  quotient $\mathbb{S}^{2n}/G$ starts at a point $(r_0, \frac{\pi}{4})$ and reaches the segment $\{r=\frac{\pi}{2}\}$ orthogonally.

\noindent Hence, by the symmetry of the solution with respect to $\{r=\frac{\pi}{2}\}$ and $\{\theta=\frac{\pi}{4}\}$, we get a smooth, compact generating curve on $\mathbb{S}^{2n}/G$. As we presented at the beginning of Section 2, the curve can generate a compact embedded CMC-hypersurface of the type $S^{n-1} \times S^{n-1} \times S^{1}$, and we prove Theorem 1.
\begin{flushright}
Q. E. D.
\end{flushright}

\section{Proof of Theorem \ref{theorem 1.2}}

\noindent In this section, we prove Theorem \ref{theorem 1.2} by modifying the analysis in Section 3.

 \noindent Let $n \ge 2$, the group $\hat G = O(n)\times O(n)\times O(n)$  acts on the sphere $\mathbb{S}^{3n-1} = \{(\vec X, \vec Y, \vec Z) : \vec X,\, \vec Y,\,\vec Z \in \mathbb{S}^{n-1} (1)\} $  in the usual way.
The quotient $\mathbb{S}^{3n-1}/\hat G  = \{ (x, y, z) \in S^{2}(1) :  x = |\vec X|=\sin(r)\cos(\theta) \geq 0,\,y = |\vec Y|= \sin(r)\sin(\theta) \geq 0, \,z=|\vec Z|= \cos(r) \geq 0\}$ is equipped with spherical coordinates $(r, \theta) \in [0, \frac{\pi}{2}] \times [0, \frac{\pi}{2}]$ and respects to the metric $g = dr^{2} + (\sin r )^{2} d\theta^{2}$.

\noindent Let $\hat C: [0, L] \longrightarrow {(\mathbb{S}^{3n-1}/\hat{G})}^{\circ},\,s \mapsto \hat{\gamma} (s)$ be a generating curve with arc length as parameter. And we take $\varphi (t_1,\dots\dots,t_{n-1}) = (\varphi_1,\dots \dots, \varphi_n)$,   $\psi (u_1,\dots\dots,u_{n-1}) = (\psi_1,\dots \dots, \psi_n)$ and $\phi (v_1,\dots\dots,v_{n-1}) = (\phi_1,\dots \dots, \phi_n)$  as three orthogonal parametrizations of the unit sphere $\mathbb{S}_{n-1}(1)$. It follows that
\begin{displaymath}
\begin{aligned}
\hat f(t_1,\dots,t_{n-1}, u_1,\dots,u_{n-1},v_1,\dots,v_{n-1}, s) =& (x(s)\varphi_1,\dots, x(s)\varphi_n, y(s)\psi_1,\dots, y(s)\psi_n,z(s)\phi_1,\dots, z(s)\phi_n),\\
\varphi_{\hat{i}}=\varphi_{\hat{i}}(t_1,\dots\dots,t_{n-1}),&\quad \varphi_1^2+ \cdots\cdots +\varphi_n^2=1,\\
\psi_{\hat{j}}=\psi_{\hat{j}}(u_1,\dots\dots,u_{n-1}),&\quad \psi_1^2+ \cdots\cdots + \psi_n^2=1,\\
\phi_{\hat{k}}=\phi_{\hat{k}}(v_1,\dots\dots,v_{n-1}),&\quad \phi_1^2+ \cdots\cdots +\phi_n^2=1, \quad 1 \leq \hat{i},\hat{j},\hat{k}, \dots \leq n
\end{aligned}
\end{displaymath}
is a parametrization of a hypersurface $\Sigma_2$ generated by the generating curve $(x(s),y(s),z(s))$, where the curve can be rewritten as $(r(s), \theta(s))$.

\noindent Furthermore, we have
\begin{displaymath}
\hat{\uppercase\expandafter{\romannumeral1}}= x^2(s)d\varphi^2 + y^2(s)d\psi^2+z^2(s)d\phi^2 + ds^2.
\end{displaymath}

\noindent And
\begin{displaymath}
\hat{\uppercase\expandafter{\romannumeral2}} = \hat{\kappa}_1\left(x^2(s)d\varphi^2\right) +\hat{\kappa}_n\left(y^2(s)d\psi^2\right)+\hat{\kappa}_{2n-1}\left(z^2(s)d\phi^2\right)+\hat{\kappa}_{3n-2}ds^2,
\end{displaymath}
where the principal curvatures
\begin{displaymath}
\begin{aligned}
\hat{\kappa}_1 = \dots = \hat{\kappa}_{n-1} =\frac{\cos r \cos \theta \sin \alpha + \sin \theta \cos \alpha }{\sin r \cos \theta},
\end{aligned}
\end{displaymath}
\begin{displaymath}
\begin{aligned}
\hat{\kappa}_{n}= \dots =\hat{\kappa}_{2n-2}  =\frac{\cos r\sin \theta\sin \alpha - \cos \theta \cos \alpha}{\sin r \sin \theta},
\end{aligned}
\end{displaymath}
\begin{displaymath}
\begin{aligned}
\hat{\kappa}_{2n-1} = \dots = \hat{\kappa}_{3n-3} =  -\tan r \sin \alpha,
\end{aligned}
\end{displaymath}
\begin{displaymath}
\begin{aligned}
\hat{\kappa}_{3n-2} =\frac{d\alpha} {ds}+ \cot r \sin \alpha.
\end{aligned}
\end{displaymath}

\noindent Then the CMC-hypersurface equation with $\hat H =\sum_{i=1}^{3n-2}\hat{\kappa}_{i} \, \equiv \, \lambda$ reduces to
\begin{equation}
\frac{d\alpha} {ds} = (2n-2)\frac{\cot (2\theta)}{\sin r} \cos \alpha + [(n-1)\tan r - (2n-1)\cot r ]\sin \alpha +  \lambda. \label{alpha1}
\end{equation}

\noindent Similarly, we need to find a compact smooth orbit for the $3 \times 3$ system
\begin{displaymath}
\left\{ \begin{array}{ll}
\frac{dr}{ds} = \cos (\alpha),  & \textrm{(\ref{r})}\\
\, \\
\frac{d\theta}{ds} = {\sin (\alpha)}/{\sin (r)},  & \textrm{(\ref{theta})}\\
\, \\
\frac{d\alpha} {ds} = (2n-2)\frac{\cot (2\theta)}{\sin r} \cos \alpha + [(n-1)\tan r - (2n-1)\cot r ]\sin \alpha +  \lambda.  & \textrm{(\ref{alpha1})}
\end{array} \right.
\end{displaymath}

\noindent {\bf{Symmetry of the solutions}}. \quad
Let $(r(s), \theta(s), \alpha(s))$ be a solution of the above system corresponding to a curve $\hat C: [0,  L] \to {(\mathbb{S}^{3n-1}/\hat G)}^{\circ}$. It is easy to verify that the solution has symmetry about $\{\theta = \frac{\pi}{4}\}$. Moreover, the symmetry of the solution with respect to $\{\tan r\cos \theta = 1\}$ (i.e.$\{x = z\}$) or $\{\tan r\sin \theta = 1\}$ (i.e.$\{y = z\}$) is also a solution of the above system by the symmetry of the coordinate components $x$, $y$, $z$.

\begin{lemma} Given any $r_0 \in (0, \arctan(\sqrt{2}))$, there exists $s_\ast > 0$ such that the system (\ref{r}), (\ref{theta}), (\ref{alpha1}) has a unique solution
\begin{displaymath}
(r, \theta, \alpha) : [0, s_\ast] \to \hat B := \{(r, \theta, \alpha) \in (0, \arctan(\sqrt{2})] \times (0, \frac{\pi}{4}] \times [-\frac{\pi}{2}, 0] :\tan r\cos \theta \leq 1 \}
\end{displaymath}
with initial data $(r, \theta, \alpha)(0) = (r_0, \frac{\pi}{4}, - \frac{\pi}{2})$ on which the solution depends continuously such that $\tan r(s_\ast)\cos \theta(s_\ast) = 1$ or $\alpha(s_\ast) = 0$. Moreover, $dr/ds \geq 0$, $d\theta/ds \leq 0$, $d\alpha/ds > 0$ and $\theta(s)> 0$ for all $s \in [0, s_\ast]$.
\end{lemma}

\noindent {\bf Proof.} Indeed, the right-hand side of equations (\ref{alpha}) and (\ref{alpha1}) for $\alpha$ have the same first and third terms and their second term is nonnegative in their domain. Therefore, we may argue exactly as in the proof of Lemma 3.1, where we simply drop the second term.
\begin{flushright}
Q. E. D.
\end{flushright}

\noindent Similarly, we condsider two cases:

\noindent {\bf Case 1}: $\lambda \in (2, \infty).$

\noindent Since $\frac{d\alpha}{ds} \geq \lambda$ still holds for all $s \in [0, s_\ast]$ and $\frac{1}{\lambda \sin(\arctan(\sqrt{2}))} < \frac{\pi}{4}$ holds for $\lambda >2$, the proofs of Lemma 3.2 and 3.3 lead to the Lemma 4.2 and 4.3 respectively.
\begin{lemma} For $n \geq 2$ and $\lambda  \in (2, \infty)$, if the initial data $r_0 \in (0, \arctan(\sqrt{2}))$ is sufficiently small, then the solution constructed in Lemma 4.1 satisfies $\alpha(s_\ast) = 0.$
\end{lemma}

\begin{lemma} For $n \geq 2$ and $\lambda  \in (2, \infty)$, let $(r, \theta, \alpha) : [0, s_\ast] \to \hat B:= \{(r, \theta, \alpha) \in (0, \arctan(\sqrt{2})] \times (0, \frac{\pi}{4}] \times [-\frac{\pi}{2}, 0] :\tan r\cos \theta \leq 1 \}$ be the solution constructed in Lemma 4.1. There exists  a constant $r_\lambda \in (0, \arctan(\sqrt{2}))$ depending only on $\lambda$ such that $\forall r_0 \in (r_\lambda, \arctan(\sqrt{2}))$, the solution constructed in Lemma 4.1 satisfies
\begin{displaymath}
0 \leq \cot(2\theta(s)) \leq \cot(\frac{\pi}{2} - \frac{2}{\lambda \sin(r_0)}), \qquad \forall s \in [0, s_\ast].
\end{displaymath}
\end{lemma}

\begin{lemma} For $n \geq 2$ and $\lambda  \in (2, \infty)$, if the initial data $r_0 \in (0, \arctan(\sqrt{2}))$ is sufficiently close to $\arctan(\sqrt{2})$, then the solution constructed in Lemma 4.1 satisfies $\tan r(s_\ast)\cos \theta(s_\ast) = 1.$
\end{lemma}

\noindent {\bf Proof.}  Since $(n-1)\tan(r(s))\sin(\alpha(s)) \leq 0, \forall s \in [0, s_\ast]$, equations (\ref{r}) and (\ref{alpha1}) imply
\begin{equation}
\begin{aligned}
\frac{d\alpha} {ds} &= \{(2n-2)\frac{\cot (2\theta)}{\sin(r)} + [(n-1)\tan(r) - (2n-1)\cot(r)]\tan(\alpha) +\frac{\lambda}{\cos(\alpha)}\}\frac{dr} {ds}  \\
&\leq [(2n-2)\frac{\cot(2\theta)}{\sin(r)} - (2n-1)\cot(r)\tan(\alpha)+\frac{\lambda}{\cos(\alpha)}]\frac{dr}{ds}.\label{34}\\
\end{aligned}
\end{equation}
With the same argument that of the proof of Lemma 3.4, we obtain $\alpha(s_\ast) < -\frac{\pi}{4}$ and  $\tan r(s_\ast)\cos \theta(s_\ast) = 1$ when $r_0$ is sufficiently close to $\arctan(\sqrt{2})$.
\begin{flushright}
Q. E. D.
\end{flushright}

\noindent {\bf Case 2}: $\lambda \in (0, 2]$.
\begin{lemma}  For $n \geq 2$ and $\lambda  \in (0, 2]$, given $r_0 \in (0, \frac{\pi}{8})$, let $(r, \theta, \alpha) : [0, s_\ast] \to  \hat B:= \{(r, \theta, \alpha) \in (0, \arctan(\sqrt{2})] \times (0, \frac{\pi}{4}] \times [-\frac{\pi}{2}, 0] :\tan r\cos \theta \leq 1 \}$ be the solution constructed in Lemma 4.1. Let $s_1 \in [0, s_\ast]$ be arbitrary, if $r(s_1) \geq 2r_0$ then
\begin{displaymath}
\theta (s_1) < \frac{\pi}{4} - \frac{1}{6n}.
\end{displaymath}
\end{lemma}

\noindent {\bf Proof.}  Since $(n-1)\tan(r(s))\sin(\alpha(s)) \leq 0, \forall s \in [0, s_\ast]$, equations (\ref{theta}) and (\ref{alpha1}) give
\begin{equation}
\begin{aligned}
\frac{d\alpha} {ds}  & = [(2n-2)\cot(2\theta)\cot(\alpha)- (2n-1)\cos(r) + \frac{\lambda\sin(r)}{\sin(\alpha)}]\frac{d\theta}{ds}+(n-1)\tan(r)\sin(\alpha)\\
\,
& \leq [(2n-2)\cot(2\theta)\cot(\alpha)- (2n-1)\cos(r) + \frac{\lambda\sin(r)}{\sin(\alpha)}]\frac{d\theta}{ds}.
\end{aligned}
\end{equation}
Thus we obtain $\theta (s_1) < \frac{\pi}{4} - \frac{1}{6n}$ with the same argument that of the proof of Lemma 3.5.
\begin{flushright}
Q. E. D.
\end{flushright}

\begin{lemma}  For $n \geq 2$ and $\lambda  \in (0, 2]$, there exists a constant $c_n>0$ depending only on n such that if the initial data $r_0 \in (0,\arctan(\sqrt{2}))$ is sufficiently small, then the solution constructed in Lemma 4.1 satisfies $\alpha(s_\ast) = 0$ and $r(s_\ast) \leq c_nr_0.$
\end{lemma}

\noindent {\bf Proof.} Let $r_0 \in (0, \frac{\pi}{8})$ be arbitrary. We may assume $r(s_\ast) > 2r_0$, otherwise the lemma directly holds. Let $0 < s_1 <s_\ast$ such that $r(s_1) = 2r_0$. From Lemma 4.5 and the monotonicity of $\theta$, we can see that
\begin{displaymath}
\theta (s) < \frac{\pi}{4} - \frac{1}{6n}
\end{displaymath}
for all $s\in[s_1, s_\ast]$. It follows that
\begin{equation}
\begin{aligned}
\frac{d\alpha} {ds} &= \{(2n-2)\frac{\cot (2\theta)}{\sin(r)} + [(n-1)\tan(r) - (2n-1)\cot(r)]\tan(\alpha) +\frac{\lambda}{\cos(\alpha)}\}\frac{dr} {ds}  \\
&\geq (2n-2)\frac{\cot(2\theta)}{\sin(r)}\frac{dr}{ds}\\
& \geq \frac{2n-2}{r}\tan(\frac{1}{3n})\frac{dr}{ds}
\end{aligned}
\end{equation}
for all $s\in[s_1, s_\ast]$.

\noindent With the same argument that of the proof of Lemma 3.6, we obtain that if the initial value $r_0>0$ is chosen such that  $c_n r_0<\frac{\pi}{4}$, then $r(s_\ast)<\frac{\pi}{4}$. And we have $\alpha(s_\ast) = 0$ by Lemma 4.1.
\begin{flushright}
Q. E. D.
\end{flushright}

\begin{lemma}  For $n \geq 2$ and $\lambda  \in (0, 2]$, if the initial data $r_0 \in (0,\arctan(\sqrt{2}))$ is sufficiently close to $\arctan(\sqrt{2})$, then the solution constructed in Lemma 4.1 satisfies $\tan r(s_\ast)\cos \theta(s_\ast) = 1.$
\end{lemma}

\noindent {\bf Proof.} Assume $r_0 \in (\arctan(\sqrt{2}- \varepsilon),\arctan(\sqrt{2}))$ for $\sqrt{2} > \varepsilon >0$. Then we have
\begin{displaymath}
\begin{aligned}
&\tan r(s)\cos \theta(s) \leq 1 \Rightarrow \cos(\theta(s)) \leq \cot(r(s)) \leq \frac{1}{\sqrt{2}-\varepsilon}\\
&\Rightarrow \theta(s) \geq \arccos(\frac{1}{\sqrt{2}-\varepsilon}) := \frac{\pi}{4} - \delta
\end{aligned}
\end{displaymath}
for all $s \in (0, s_\ast)$. Hence,
\begin{displaymath}
\theta(s) \geq \frac{\pi}{4} - \delta,\quad 0 \leq \cot(2\theta) \leq  \tan(2\delta) := b, \quad \cot(\alpha) \leq 0,\quad \sin(\alpha)<0, \quad \cos(r) \in  (0, 1), \quad \frac{d\theta}{ds} \leq 0,
\end{displaymath}
for all $s \in [0, s_\ast)$. Since
\begin{displaymath}
\begin{aligned}
\frac{d\alpha} {ds}  & = [(2n-2)\cot(2\theta)\cot(\alpha)- (2n-1)\cos(r) + \frac{\lambda\sin(r)}{\sin(\alpha)}]\frac{d\theta}{ds}+(n-1)\tan(r)\sin(\alpha)\\
\,
& \leq [(2n-2)\cot(2\theta)\cot(\alpha)- (2n-1)\cos(r) + \frac{\lambda\sin(r)}{\sin(\alpha)}]\frac{d\theta}{ds}\\
& \leq (2n-1)(b\cot(\alpha) - 1 + \frac{1}{\sin(\alpha)})\frac{d\theta}{ds},
\end{aligned}
\end{displaymath}
choosing $s_1 = s_\ast$ and $\varepsilon \in (0,\sqrt{2}-\frac{1}{\cos(\frac{\pi}{4}-\frac{1}{6n})})$, we get $\delta \leq \frac{1}{6n}$. With the same argument that of the proof of Lemma 3.5, we obtain  $\alpha (s_\ast) \leq (2n-1)(b^2+1)\delta-b\ln(\frac{b+1}{2})-\frac{\pi}{2}+\frac{1}{b}\ln(1+b)$, then we have $\alpha(s_\ast) < -\frac{1}{2}$  with $\varepsilon$ is suffeicently small. Moreover, since $f(y)= \ln(\frac{y+1}{(1-b)y+1+b})$ is increasing in $(0,1)$, let $y_1 = \tan(\frac{1}{4})$, it follows that
\begin{displaymath}
\begin{aligned}
\alpha (s_\ast) &\leq (2n-1)(b^2+1)\delta-b\ln(\frac{b+1}{2})-\frac{\pi}{2}- \frac{1}{b}f(-\tan(\frac{\alpha(s_\ast)}{2}))\\
&< (2n-1)(b^2+1)\delta-b\ln(\frac{b+1}{2})-\frac{\pi}{2}+\frac{1}{b}\ln(1+\frac{1-y_1}{1+y_1}b)\\
&< -\frac{\pi}{4}.
\end{aligned}
\end{displaymath}

\noindent Thus we have $\alpha(s_\ast) < -\frac{\pi}{4}$ and $\tan r(s_\ast)\cos \theta(s_\ast) = 1$ for $r_0 \in (\arctan(\sqrt{2}- \varepsilon),\arctan(\sqrt{2}))$  when $\varepsilon$ is suffeicently small.
\begin{flushright}
Q. E. D.
\end{flushright}

\noindent {\bf Proof of Theorem Theorem \ref{theorem 1.2}}. It is easy to prove that
\begin{displaymath}
\beta(r,\theta) = -\arctan(\sin(r)\sin(\theta))
\end{displaymath}
is the angle between the normal vector $(1, -\sin \theta)$  of the curve $\{ \tan r\cos \theta = 1\}$ and the vector $\partial/\partial r$ at the point $(r,\theta) = (r,\arccos(\cot(r)))$.

\noindent And we have found that for $n \geq 2$ and $\lambda  > 0$, there exist $r_0',\,r_0'' \in (0, \arctan(\sqrt{2}))$ with $r_0' < r_0''$ such that the trajectories emanating from  $(r_0', \frac{\pi}{4}, -\frac{\pi}{2})$ and $(r_0'', \frac{\pi}{4}, -\frac{\pi}{2})$ leave the domain $\hat B := \{(r, \theta, \alpha) \in (0, \arctan(\sqrt{2})] \times (0, \frac{\pi}{4}] \times [-\frac{\pi}{2}, 0] :\tan r\cos \theta \leq 1 \}$ from $\{\alpha = 0\}$ and from $\{\tan r\cos \theta = 1\} \cap \{\alpha < -\frac{\pi}{4}\}$ respectively.
Similar to the proof of Theorem 1, we can prove that there exists a $r_0''' \in (r_0', r_0'')$ (by the continuity of the system of ODEs), for which, the trajectory emanating from $(r_0''', \frac{\pi}{4}, -\frac{\pi}{2})$ leaves the domain $\hat B$ from $\{\tan r\cos \theta = 1\}  \cap  \{\alpha= 0\}$.

\noindent Since $\beta(r,\theta) \in [-\frac{\pi}{4},0]$, we can find a trajectory that exits the domain $\hat B$ from $\{\tan r\cos \theta = 1\}  \cap  \{\alpha= \beta(r,\theta)\}$. The projection of this trajectory on the quotient ${(\mathbb{S}^{3n-1}/\hat G)}^{\circ}$ starts at a point $(r_0, \frac{\pi}{4})$ and reaches the segment $\{\tan r\cos \theta =1 \}$ orthogonally.

\noindent By the symmetry of the solution $(r(s), \theta(s), \alpha(s))$, we get a smooth, compact generating curve, which implies the desired conclusion.
\begin{flushright}
Q. E. D.

\end{flushright}

\end{document}